\documentclass[12pt]{amsart}
\usepackage{amsmath}
\usepackage{amssymb}

 \newtheorem{thm}{Theorem}[section]
 \newtheorem{cor}[thm]{Corollary}
 \newtheorem{lem}[thm]{Lemma}
 \newtheorem{prop}[thm]{Proposition}
 \theoremstyle{definition}
 \newtheorem{defn}[thm]{Definition}
 \theoremstyle{remark}
 \newtheorem{rem}[thm]{Remark}
 \numberwithin{equation}{subsection}
 
 \DeclareMathOperator{\Gcd}{gcd}
 \DeclareMathOperator{\tr}{tr}
 \DeclareMathOperator{\Diag}{Diag}
 \DeclareMathOperator{\Gal}{Gal}
 \DeclareMathOperator{\N}{N}
 
 \DeclareMathOperator{\C}{C}
 \DeclareMathOperator{\h}{H}
 \DeclareMathOperator{\Char}{Char}
 \DeclareMathOperator{\Br}{Br}

\begin{document}





\title[All dihedral division algebras of degree five are cyclic]
{\textbf{All dihedral division algebras of degree five are cyclic}}

\author{ Eliyahu Matzri }

\address{Department of Mathematics, Bar-Ilan University, Ramat-Gan, 52900, Israel}

\email{elimatzri@gmail.com}

\thanks{The author thanks his supervisors,
        L.H. Rowen and U. Vishne, for many interesting and motivating
        talks and for supporting this work through BSF grant no.~2004-083.}







\begin{abstract}
In \cite{AAB95} Rowen and Saltman proved that every division algebra
which is split by a dihedral extension of degree $2n$ of the center,
$n$ odd, is in fact cyclic. The proof requires roots of unity of
order $n$ in the center. We show that for $n=5$, this assumption can
be removed. It then follows that ${}_{5\!\!\!\:}\Br(F)$, the
$5$-torsion part of the Brauer group, is generated by cyclic
algebras, generalizing a result of Merkurjev \cite{AAC95} on the $2$
and $3$ torsion parts.
\end{abstract}

\maketitle


\section{\textbf{Mathematical background}}

\bigskip
We begin with basic notions needed for this work and refer the
reader to \cite{RRB95} or \cite{NN59} for more details.\\
Let $R$ be a ring and let $\C(R)=\{r\in R \mid rx=xr \hbox{ \ }
\forall x\in R \}$ denote the center of $R$.
\begin{defn}
A ring $R$ will be called a simple ring if $R$ has no non-trivial
two-sided ideals. In particular $R$ is a division ring if every
nonzero element is invertible.
\end{defn}
\smallskip

\begin{rem}
Notice that if $R$ is simple, its center is naturally a field.
\end{rem}

\smallskip
\begin{defn}
An $F$-algebra $R$ is called an $F$-central simple algebra if $R$ is
simple with $\C(R)=F$ and $\dim_F(R)<\infty$.
\end{defn}
\begin{rem}
Every $F$-central simple algebra $A$ has $\dim_F(A)=n^2$, and we
define the degree of $A$, denoted $\deg(A)$, to be $n$.
\end{rem}

By Wedderburn's Theorem every $F$-central simple algebra is of the
form $M_n(D)$, where $D$ is a division algebra with center $F$. \ \\
\smallskip
The Brauer group of a field $F$, denoted $\Br(F)$, is the set of
isomorphism classes of $F$-central simple algebras modulo the
following relation: two central simple algebras $A,B$ are equivalent
if and only if there exist natural numbers $n,m$ such that
$M_n(A)\cong M_m(B)$.
\smallskip

\begin{prop}
Let $D$ be an $F$-central division algebra of degree $n$, and $K$ a
subfield of $D$, then $K$ is a maximal subfield if and only if
$[K:F]=n$.
\end{prop}
\begin{defn}
A crossed product is an $F$-central simple algebra $A$ of degree $n$
containing a commutative $F$-subalgebra $C$ Galois over $F$, with
$[C:F]=n$. Note that if $A$ is a division algebra then $C$ is a
maximal subfield of $A$.
\end{defn}
\smallskip
\begin{defn}
Let $D$ be an $F$-central division algebra of degree $n$. We will
say that $D$ is split by a group $G$ if $D$ contains a maximal
subfield $K$ with Galois closure $E$ such that $\Gal(E/F)=G$.
\end{defn}
\smallskip
\begin{thm}
Let $A$ be a crossed product where $K\subset A$ is a maximal
subfield with Galois group $\Gal(K/F)=G$. Then $A$ has the following
description: $A=\mathop \oplus \limits_{\sigma  \in G}Kx_{\sigma} $
as a left $K$-vector space, and multiplication in $A$ is according
to the rules:
$$x_\sigma k=\sigma(k)x_\sigma \hbox{ \ } \forall k\in K$$ and
$$x_\sigma x_\tau=c(\sigma,\tau)x_\tau x_\sigma$$ where $c\in \h
^2(G,K^{\times})$ is a $2$-cocycle. In this case $A$ is denoted
$A=(K,G,c)$.

\end{thm}
\smallskip
\begin{rem}
If $G=\left<\sigma\right>$ we can give a simpler representation of
$A$ as follows:\\ $A=\mathop  \oplus \limits_{i = 0}^{n - 1}Kx^i$ as
a left $K$-vector space, where $n=\deg(A)=\left | G\right |$ and the
multiplication is according to the rules:
$$xk=\sigma(k)x \hbox{ \ } \forall
k\in K$$ and $$x^ix^j=\left\{
                        \begin{array}{ll}
                          x^{i+j}, & i+j<n \\
                          \beta x^{i+j-n}, & i+j \geq n
                        \end{array}
                      \right.
$$ \\
In this case, $A$ is denoted as $A=(K,\sigma,\beta)$.
\end{rem}
\smallskip
\begin{rem}
If $F$ contains a primitive $n$-th root of unity $\rho$, we can give
an even simpler description of $A$ (since then $K=F[x\mid
x^n=\alpha\in F]$) as follows: $$A=F[x,y\mid
x^n=\alpha;y^n=\beta;xy=\rho_nyx]\hbox{ \ \ \ } \alpha,\beta\in F$$
\end{rem}
\smallskip

\section{\textbf{Some preliminary results}}
\smallskip

In this section we briefly repeat the arguments of Rowen and Saltman
in \cite{AAB95} but we do not assume $F$ contains roots of unity.

The situation we will be handling is the following: \\ $D/F$ is a
central simple algebra of odd degree $n$ having a maximal subfield
$K\subset D$ with Galois closure $E\supset K\supset F$, such that
$$\Gal(E/F)=D_n=\left<\sigma, \tau : \sigma^n=\tau^2=1 \ ,\tau \sigma
\tau = \sigma^{-1}\right>,$$ and $K=E^{\left<\tau\right>}$.

Extending scalars to $E^{\left<\sigma\right>}$, we may view
$E\subset D'=D\otimes E^{\left<\sigma\right>}$. Now
$\Gal(E/E^{\left<\sigma\right>})=\left<\sigma\right>$, i.e. $D'$
is cyclic, so we have an element $\beta\in D'$ such that\\
$(1) \hbox{ \ \ \ \ \ \ \ \ \ \ \ \ \ \ \ \ \ \ \ \ \ }\beta^{-1}x\beta=\sigma (x)\hbox{ \ \ \ } \forall x\in E.$\\
In particular $\beta^n \in E^{\left<\sigma\right>}$. Notice that
$\tau$ can be extended to $D'=D\otimes E^{\left<\sigma\right>}$ by
its action on $E^{\left<\sigma\right>}$, that is, we write $\tau$
instead of $1\otimes\tau$.

\smallskip

\begin{lem}
We may assume that $\tau(\beta)=\beta^{-1}$.
\end{lem}

\smallskip

\begin{proof}
Applying $\tau$ to $(1)$ yields
$$\tau(\beta)^{-1}\tau(x)\tau(\beta)=\sigma^{-1}(\tau(x)) \hbox{ , \ \ \ } \forall x\in E.$$
Now since $\tau$ is an automorphism of $E$, $\tau(x)$ runs over all
elements of $E$, and thus
$$\tau(\beta)^{-1}y\tau(\beta)=\sigma^{-1}(y) \hbox{ , \ \ \ }
\forall y\in E$$ that is $\tau(\beta)$ acts on $E$ as $\sigma^{-1}$.
Now define $\beta '=\beta^r \tau(\beta)^{-r},$ where $r=(n+1)/2,$
and compute that $\tau(\beta')=\beta'^{-1}$, and $\beta'$ acts on
$E$ as $\sigma$.
\end{proof}

\smallskip

Let $P_t(X)=X^n+\sum_{i=1} ^{n} c_i(t)X^{n-i}$ denote the
characteristic polynomial of $t\in D'$. Note that $c_1(t)=-\tr(t)$
and $c_n(t)=(-1)^n\N(t)$ where $tr(t)$ and $N(t)$ are the reduced
trace and norm of $t$.

\smallskip

\begin{lem}\label{t}
Let $t=\beta^i e,$ for $e\in E$ and $0<i<n$, $i\neq 0$. Then
$tr(t)=0$.
\end{lem}

\smallskip

\begin{proof}
Let $d=\Gcd(i,n)$.\\
Clearly we have $t^{n/d}=\beta^{ni/d}\N_{\sigma^i}(e)\in
E^{\left<\sigma^i\right>}$ where $\N_{\sigma^i}$ is the norm from
$E$ to $E^{\left<\sigma^i\right>}$. Now
$[E:E^{\left<\sigma^i\right>}]=n/d,$ implying
$P(X)=X^{n/d}-\beta^{ni/d}\N_{\sigma^i}(e)$ is the characteristic
polynomial of $t$, hence $\tr_{E/E^{\left<\sigma^i\right>}}(t)=0$
which implies $\tr_{E/F}(t)=0$.
\end{proof}

\smallskip

\begin{lem}
Let $t=(\beta+\beta^{-1})e$ for $e\in E$. Then the coefficients of
$P_t(X)$ satisfy $c_i(t)=0$ for every odd $0<i<n$.
\end{lem}

\smallskip

\begin{proof}
Notice that for $i$ odd, $t^i$  is a sum of elements of the form
$a\beta^s$ where $a\in E$ and $s$ odd, $-n< s < n$, so by \ref{t}
and Newton's identities we are done in the characteristic zero case.
For the general case, we refer the reader to \cite{AAB95} where the
main idea is that you can form a model for this situation in the
form of an Azumaya algebra and then use a specialization argument.
\end{proof}

\smallskip

\begin{cor}\label{b}
There is an element $t\in D$ such that for every $e\in E$ (and so
also for $k\in K\subset E$), $c_i=0$ for every odd $0<i<n$ in
$P_{te}(X)$.
\end{cor}

\smallskip

\begin{proof}
Since $D=D'^{\left<\tau\right>}$ we have $t=\beta+\beta^{-1}$ is the
desired element.
\end{proof}

\smallskip

\begin{rem}\label{char}
Notice that if $n=p$ is prime $\Char(F)=p$, the element
$t=\beta+\beta^{-1}\in D$ we found satisfies $t^p\in F$ and $t\notin
F$ and so by a theorem of Albert in the ``special results'' chapter
of his seminal book \cite{dB22}, which is knows as Albert's
cyclicity criterion, $D$ is cyclic (this is not a new result as J.P
Tignol and P. Mammone did this for any field $F$ with $\Char(F)\mid
n$ in \cite{Dev68} using the corestriction, but it shows that the
proof of Rowen and Saltman also applies to this case).
\end{rem}

\bigskip
\section{\textbf{The case $n=5$ }}
\smallskip

Now we would like to focus on the particular case where $n=5$. The
main tool we will be using is the following proposition taken from
[$3$, Proposition $2.2$].
\smallskip

\begin{prop}\label{deg3}
Let $G(x_1,...,x_n)$ be a homogeneous form of degree $3$ defined
over a field $F$. If $G$ has a solution, $\alpha \in K^{(n)}$,
defined over a quadratic extension $K$ of $F$, then $G$ has a
solution defined over $F$.
\end{prop}
\smallskip
\begin{proof}
The proof in \cite{dB93} uses basic intersection theory which we
will not use, instead we will give an algebraic proof (which is
actually a translation of the proof in \cite{dB93}) which will
enable us to find an explicit solution in section $3$. Since
$[K:F]=2$ the solution $\alpha$ has the following form:
$\alpha=(\alpha_1+\beta_1t,...,\alpha_n+\beta_nt)$ where
$\alpha_i,\beta_i\in k$, and $t\in K$ such that $K=F[t]$. Now
specialize $G(x_1,...,x_n)$ to
$G(\alpha_1+\beta_1Z,...,\alpha_n+\beta_nZ)$, denoting it by $g(Z)$.
Notice that the coefficient of $Z^3$ in $g(Z)$ is
$G(\beta_1,...,\beta_n)$ hence if $G(\beta_1,...,\beta_n)=0$ we have
a solution defined over $F$ else $g(Z)$ is a degree $3$ polynomial
defined over $F$. Since $g(t)=0$ we get that $g(Z)=cm_t(Z)(Z-w)$,
where $c=G(\beta_1,...,\beta_n)$ and $m_t(Z)$ is the minimal
polynomial of $t$ over $F$. Now $c$, $g(Z)$ and $m_t(Z)$ are defined
over $F$ hence $w$ is in $F$ and clearly
$G(\alpha_1+\beta_1w,...,\alpha_n+\beta_nw)=g(w)=0$ so we have found
a solution $\gamma=(\alpha_1+\beta_1w,...,\alpha_n+\beta_nw)\in
F^n$.
\end{proof}

\smallskip

\begin{thm}\label{yyy}
Let $D$ be a division algebra of degree $5$ split by the group $D_5$
then $D$ is cyclic.
\end{thm}

\smallskip

\begin{proof}
In view of remark \ref{char}, we may assume $\Char(F)\neq5$. First
we remark that by Albert's cyclicity criterion it is enough to find
an element $t\in D-F$ such that $t^5\in F$, that is $c_i=0$ for
every $0<i<n$. Now by \ref{b} we have $t\in D$ with the property
$c_i(te)=0$ for every odd $0<i<n$ and $\forall e\in E$. Now since
$P_{t^{-1}}(x)=-N(t)^{-1}P_t(x^{-1})x^5$ we have $c_i(et^{-1})=0$
for every even $0<i<n$ and $\forall e\in E$. Hence we are left with
finding a solution for $c_1(et^{-1})=0$ (which is linear) and
$c_3(et^{-1})=0$ (which is cubic) in the five dimensional vector
space $Et^{-1}$. Define $V:=\{et^{-1}\in Et^{-1} \mid
c_1(et^{-1})=0\}$, which is a four dimensional subspace of
$Et^{-1}$. We have to find a solution for $c_3(v)=0$ in $V$. Let us
add a fifth root of unity to $F$, which is either a quadratic
extension or a chain of two quadratic extensions. After this
extension we are in the case of Rowen and Saltman where they gave an
explicit element whose fifth power is in $F$ which was
$(v+v^{-1})t^{-1}$, where $v\in E$. This element is clearly in
$V\otimes_F F[\rho_5]$. Now by \ref{deg3} since $c_3(v)$ is
homogeneous of degree $3$, we have a solution after either one or
two quadratic extensions. Thus, we have a solution before the
extension and we are done.
\end{proof}

\smallskip

\begin{rem}
If the fifth root of unity is in a quadratic extension of $F$, we
know $D$ is cyclic by a theorem of Vishne [$10$, Theorem $13.6$] and
D. Haile, M. A. Knus, M. Rost, J. P. Tignol \cite{Gam90}, so what
actually is new is the last case of $[F[\rho]:F]=4$.
\end{rem}

\bigskip

\section{\textbf{A generic example}}
Fixing $p$ let $K=F[\rho_p]$ and denote
$\Gal(K/F)=\left<\tau\right>$. In [$6$, Theorem $2$] Merkurjev
proves that ${}_{p\!\!\!\:}\Br(F)$ is generated by $F$-central
simple algebras, $A$, of degree $p$ such that $A\otimes K \simeq
(\alpha,\beta)$
where $K[\sqrt[p]{\alpha}]$ is cyclic over $K$ Galois over $F$.\\
In \cite{Enf87} Vishne calls these algebras quasi-symbols and gives
more details about them including generic examples. We will show
that for $p=5$ these algebras are cyclic and conclude that
${}_{5\!\!\!\:}\Br(F)$ is generated by cyclic algebras.


\subsection{\textbf{A generic Quasi-symbol of degree $5$ }} \ \\
\smallskip
For $p=5$ we have two possibilities for $[K:F]$. The first is
$[K:F]=2$; in this case Vishne shows that every quasi-symbol is
cyclic.  The second case is $[K:F]=4$; in this case every
quasi-symbol $A$ has one of the following forms (after extending
scalars to $K$):
\begin{enumerate}
  \item $A\otimes K=(\alpha,\beta)$, where $\alpha\in F$ and
  $\tau(\beta)\equiv\beta^2 \hbox{ \ }\pmod{K^{\times^5} }$.
  \item $A\otimes K=(\alpha,\beta)$, where  $\tau(\alpha)=\alpha^{-1}$ and
  $\tau(\beta)\equiv\beta^{-2} \hbox{ \ }\pmod{K^{\times^5}}$.
\end{enumerate}

The first kind is known to be cyclic by [10, Theorem 10.3]. So we
are left with the second kind for which Vishne gives the following
generic construction which we will show is cyclic. Thus every
quasi-symbol of degree $5$ is cyclic and hence, by [$6$, Theorem
$2$] we conclude that ${}_{5\!\!\!\:}\Br(F)$ is generated by cyclic
algebras.

\bigskip

Let $k_0$ be a field of characteristic $\neq 5$ and $k=k_0[\rho]$
where $\rho$ is a fixed primitive fifth root of unity,
$\Gal(k/k_0)=\left<\tau\right>$ where $\tau(\rho)=\rho^2$. Set
$K=k(a,b,\eta)$ a transcendental extension and extend $\tau$ to $K$
by
$$\tau(a)=a^{-1}, \hbox{ \ \ \ \ \ } \tau(b)=\eta^5b^{-2},  \hbox{ \ \ \ \ \
}\tau(\eta)=\eta^2b^{-1}.$$ Notice that we still have $\tau^5=1$.
Define $F=K^{\left<\tau\right>}$ and
$$D=(a,b)_K=K[x,y\mid x^5=a,\hbox{ \ \ } y^5=b,\hbox{ \ \ }
yxy^{-1}=\rho x],$$ and extend $\tau$ to $D$ by
$\tau(x)=x^{-1},\hbox{ \ \ }\tau(y)=\eta y^{-2}$.
Notice that $\tau^2(\eta)=\eta^{-1}$ and $\tau^2(y)=y^{-1}$.\\
Now define $D_0=D^{\left<\tau\right>}$; $D_0/F$ is the generic
quasi-symbol of degree $5$ of the second type.

\begin{rem}
Vishne's construction is much more general and we specialized it to
the above case, for the general construction we refer the reader to
\cite{Enf87}.
\end{rem}

\smallskip

\begin{prop}
$D_0$ is split by $D_5$.
\end{prop}
\smallskip
\begin{proof}
Notice that $\Gal (K[y]/F)=C_5 \rtimes C_4=\left<\sigma\right>
\rtimes\left<\tau\right>$ and now we will see how $\tau$ acts on
$\sigma$. Applying $\tau$ to $x^{-1}tx=\sigma(t)$, which holds for
every $t\in K[y]$, yields
$\tau(\sigma(t))=\tau(x^{-1})\tau(t)\tau(x)=x\tau(t)x^{-1}=\sigma^{-1}(\tau(t))$
and so we get $\tau\sigma\tau^{-1}=\sigma^{-1}$. Hence $\tau^2$ is a
central element in $\Gal (K[y]/F)$  and it is clear that
$E=K[y]^{\left<\tau^2\right>}\subset K[y]$ is Galois over $F$ with
$\Gal(E/F)=D_5=\left<\sigma\right> \rtimes\left<\tau\right>$ and we
are done.
\end{proof}

\smallskip
\begin{cor}
$D_0$ is cyclic.
\end{cor}
\smallskip

In \cite{AAC95} Merkurjev proves the following theorem:
\begin{thm}
Let $F$ be a field. ${}_{n\!\!\!\:}\Br(F)$ is generated by cyclic
algebras, for $n=2,3$.
\end{thm}

\smallskip

Now as a result of the above we can extend Merkurjev's theorem to
$n=5$ and get

\begin{thm}
${}_{5\!\!\!\:}\Br(F)$ is generated by cyclic algebras.
\end{thm}
\smallskip
\begin{proof}
By section $8$ of \cite{Enf87} ${}_{5\!\!\!\:}\Br(F)$ is generated
by quasi-symbols of degree $5$, and so we are done.
\end{proof}

\bigskip
\subsection{\textbf{Finding an explicit solution}} \ \\
\smallskip
Since the above example is a generic one, it would be nice to give
an explicit element with fifth power in $F$, which is what we do now
by going over the general proof.

Let $P_t(X)=X^n+\sum_{i=1} ^{n} c_iX^{n-i}$ denote the
characteristic polynomial of $t\in D_0$.\\
$V=(x+x^{-1})^{-1}K[y]^{\left<\tau\right>}$ is a $5$-dimensional
$F$-subspace of $D_0$, satisfying $c_2(v)=c_4(v)=0$ for all $v\in
V$; and we want to find a solution in $V$ for
$tr(Z)=c_1((x+x^{-1})^{-1}Z)=0$ and $G(Z)=c_3((x+x^{-1})^{-1}Z)=0$.
Extending scalars from $F$ to $F[\rho+\rho^{-1}]$, we have the
solutions $Z_1=y+y^{-1}=\alpha+\beta(\rho+\rho^{-1})$ and
$Z_2=\tau(Z_1)=\alpha+\beta\tau(\rho+\rho^{-1})=\alpha+\beta\tau(\rho^2+\rho^{-2})$
where $\alpha=(\alpha_1,...,\alpha_5), \beta
=(\beta_1,...,\beta_5)\in K[y]^{\left<\tau\right>}$ so
$\alpha_i,\beta_i\in F$. Now define the following line:
$L=\{\alpha+\beta t\}=\{(\alpha_1+\beta_1t,...,\alpha_5+\beta_5t)\}$
defined over $F$.
\smallskip
\begin{prop}
For every $l\in L$ we have $tr(l)=0$.
\end{prop}
\smallskip
\begin{proof}
By standard linear algebra, $L\cap \{tr(Z)=0\}$ is either one point
or the whole line $L$; since $Z_1,Z_2\in L\cap \{tr(Z)=0\}$, we get
$L\cap \{tr(Z)=0\}=L$ and we are done.
\end{proof}

\smallskip

Now let us study the variety $\{G(Z)=0\}\cap L$. First we need to
compute $G(Z)$. In order to do that we use the representation of $D$ induced by right multiplication on\\
$D=K[y]+K[y]x+K[y]x^2+K[y]x^3+K[y]x^4$, namely
$$x\longrightarrow \left(
                     \begin{array}{ccccc}
                       0 & 0 & 0 & 0 & a \\
                       1 & 0 & 0 & 0 & 0 \\
                       0 & 1 & 0 & 0 & 0 \\
                       0 & 0 & 1 & 0 & 0 \\
                       0 & 0 & 0 & 1 & 0 \\
                     \end{array}
                   \right)$$
$$m\in K[y]\longrightarrow
\Diag(m,\sigma(m),\sigma^2(m),\sigma^3(m),\sigma^4(m))$$

Now the minimal polynomial of $x+x^{-1}$ is
$$\lambda^5-5\lambda^3+5\lambda-(a+a^{-1})$$ hence
$$(x+x^{-1})^{-1}=((x+x^{-1})^4-5(x+x^{-1})^2+5)(a+a^{-1})^{-1}=(a+a^{-1})^{-1}(x^4+x^{-4}-x^2-x^{-2}+1)$$
implying $$(x+x^{-1})^{-1}\longrightarrow (a+a^{-1})^{-1} \left(
                                        \begin{array}{ccccc}
                                          1 & a & -1 & -a & 1 \\
                                          a^{-1} & 1 & a & -1 & -a \\
                                          -1 & a^{-1} & 1 & a & -1 \\
                                          -a^{-1} & -1 & a^{-1} & 1 & a \\
                                          1 & -a^{-1} & -1 & a^{-1} & 1 \\
                                        \end{array}
                                      \right)$$

Now when we compute the characteristic polynomial of
$(x+x^{-1})^{-1}m$ we get that:\\
 \ \\

$$c_3((x+x^{-1})^{-1}m)=(a+a^{-1})^{-1}(m\sigma(m)\sigma^2(m)+\sigma(m)\sigma^2(m)\sigma^3(m)
+\sigma^2(m)\sigma^3(m)\sigma^4(m)+$$
$\sigma^3(m)\sigma^4(m)m+\sigma^4(m)m\sigma(m))=(a+a^{-1})^{-1}\tr_\sigma(m\sigma(m)\sigma^2(m))
$. \\
 \ \\
Yielding $F(Z)=(a+a^{-1})^{-1}\tr_\sigma(Z\sigma(Z)\sigma^2(Z))$
\\Now clearly $\{F(Z)=0\}\cap L$ is defined over $F$ by the
polynomial
\\ $f(t)=F(\alpha+\beta t)=(a+a^{-1})^{-1}\tr_\sigma (\alpha+\beta
t)\sigma(\alpha+\beta t)\sigma^2(\alpha+\beta t))=
(a+a^{-1})^{-1}\tr_\sigma(\beta\sigma(\beta)\sigma^2(\beta)t^3+...)=F(\beta)t^3+...$\\
But we know two solutions for $f(t)$, namely $t_1=\rho+\rho^{-1}$
and $t_2=\rho^2+\rho^{-2}$, so we get
$f(t)=F(\beta)(t-t_1)(t-t_2)(t-t_3)$. Now since $f(t)$ and
$F(\beta)(t-t_1)(t-t_2)$ are defined over $F$, we get $t_3\in F$.\\
Explicitly $f(0)=-t_1t_2t_3F(\beta)$ implies
$t_3=\frac{-f(0)}{t_1t_2F(\beta)}=\frac{f(0)}{F(\beta)}=\frac{F(\alpha)}{F(\beta)}$
is in $F$. Hence we get:\\
\begin{thm}
The element $w=(x+x^{-1})^{-1}(\alpha+\beta
\frac{F(\alpha)}{F(\beta)})\in D_0-F$ satisfies $w^5\in F$.
\end{thm}
 \ \\
Now we are left with solving for $\alpha,\beta$ from the two
equations:
$$y+y^{-1}=\alpha+\beta(\rho+\rho^{-1})$$
$$\eta y^{-2}+\eta^{-1} y^2=\tau(y+y^{-1})=\alpha+\beta(\rho^2+\rho^{-2})$$
Hence \\
$$\beta=\frac{y+y^{-1}-\eta y^{-2}-\eta^{-1}
y^2}{\rho+\rho^{-1}-\rho^2-\rho^{-2}}$$
$$\alpha=y+y^{-1}-\beta(\rho+\rho^{-1})$$

 \ \\
 \ \\
  \ \\

\subsection{\textbf{The general case}} \ \\
\smallskip
We will now show that the above solution for the case of
quasi-symbols, where we do decent from $F[\rho+\rho^{-1}]$ to $F$ is
valid for the general case of $D_5=\left<\sigma, \tau :
\sigma^5=\tau^2=1,\tau \sigma \tau = \sigma^{-1}\right>$ division
algebras, where we need to decent from $F[\rho]\otimes E^{\left<
\sigma \right>} $ to $F$. The situation is the following: we look
for a solution to $c_3(t)=c_1(t)=0$ where $c_i(t)$ are as in section
$3$ and $t\in (\beta+\beta^{-1})^{-1}E^{\left< \tau \right>}$. Let
$\Gal(F[\rho]/F)=\left< \pi \right>$; hence $\Gal (E\otimes
F[\rho]/F)=D_5\times \left< \pi \right>$  and so after extending
scalars to $F[\rho]$ we want a solution in
$(\beta+\beta^{-1})^{-1}(E\otimes F[\rho])^{\left< \tau
\right>\times\left< \pi \right>}$, which will then be defined over
$F$.
\smallskip
\begin{prop}
We may assume $v+v^{-1}\in (E\otimes F[\rho])^{\left< \tau
\right>\times\left< \pi^2 \right>}$, for $v$  as in the proof of
theorem \ref{yyy}.
\end{prop}
\smallskip
\begin{proof}
Since $v=x^r\tau(x)^{-r}$, where $x$ is any eigenvector of $\sigma$
with eigenvalue $\rho$, we may write $x=\sum_{i=0}^4
\rho^{-i}\sigma^i(k)$ for $k\in E^{\left< \tau \right>\times \left<
\pi \right>}$. Now $\tau(x)=\pi^2(x)$ and so
$\tau(v)=\tau(x)^rx^{-r}=\pi^2(x)^rx^{-r}=\pi^2(x^r\pi^2(x)^{-r})=\pi^2(x^r\tau(x)^{-r})=\pi^2(v)$
implying $\tau(v+v^{-1})=v+v^{-1}$, hence $v+v^{-1}$ is in
$(E\otimes F[\rho])^{\left< \tau \right>\times\left< \pi^2
\right>}$, as desired.
\end{proof}
\smallskip
Now it is clear that after extending scalars to $F[\rho+\rho^{-1}]$
we have the solution $(\beta+\beta^{-1})^{-1}(v+v^{-1})$ and so we
are in the same situation as in the quasi-symbol case, hence the
above solution is valid for the general case too .


\end{document}